\documentclass[letterpaper, 10pt]{article}
\usepackage{graphicx}
\usepackage{float}
\usepackage[inner=0.75in, outer=0.75in, top=0.75in, bottom=1in, paperheight=11in, paperwidth=8.5in]{geometry}

\usepackage{array}
\usepackage{rotating}
\usepackage{multirow}

\usepackage{mathtools}
\usepackage{shuffle}
\usepackage[table, svgnames, dvipsnames,xcdraw]{xcolor}
\usepackage[square, numbers, comma, sort&compress]{natbib}
\usepackage{verbatim}
\usepackage{bm}
\usepackage{amsmath}
\usepackage{amsfonts}
\usepackage{amssymb}
\usepackage{amsthm}

\usepackage{epstopdf}
\usepackage[utf8]{inputenc}
\usepackage[english]{babel}
\usepackage{color}
\usepackage{forest}

\usepackage{pdfpages}
\usepackage{pgfplots}
\usepackage{pgfplotstable}
\usepgfplotslibrary{patchplots}
\pgfplotsset{compat=1.15}

\usepackage{microtype}
\usepackage{url}
\usepackage{subfiles}
\usepackage{caption}
\usepackage{subcaption}
\usepackage{parskip}
\usepackage[section]{placeins}

\usepackage{tikz}
\usepackage{tikz-cd}
\usetikzlibrary{arrows,calc,shapes,decorations.pathreplacing}
\usetikzlibrary{decorations.markings}
\usetikzlibrary{trees,positioning,fit,shapes}
\usetikzlibrary{backgrounds}
\usetikzlibrary{patterns, positioning, arrows}

\usepackage{framed}

\colorlet{shadecolor}{blue!20}

\usepackage[obeyFinal,textwidth=0.55in]{todonotes}
\usepackage[hidelinks]{hyperref}
\usepackage[capitalize]{cleveref}

\usepackage{enumitem}

\theoremstyle{plain}
\newtheorem{thm}{Theorem}[section]
\theoremstyle{definition}
\newtheorem{defn}[thm]{Definition}

\newtheorem{lemm}[thm]{Lemma}
\newtheorem{corr}[thm]{Corollary}

\makeatletter
\newcommand\setveclength[3]{
  \pgfpointdiff{\pgfpointanchor{#2}{center}}{\pgfpointanchor{#3}{center}}
  \pgfmathveclen{\pgf@x}{\pgf@y}
  \edef#1{\pgfmathresult}
}
\makeatother

\begin{document}	

\title{\bf On topological data analysis for SHM; an introduction to persistent homology}	
\author{T.\ Gowdridge, N.\ Dervilis, K.\ Worden \\
        Dynamics Research Group, Department of Mechanical Engineering, University of Sheffield \\
        Mappin Street, Sheffield S1 3JD, UK
	   }
	\date{}
    \maketitle
	\thispagestyle{empty}

\section*{Abstract}
This paper aims to discuss a method of quantifying the 'shape' of data, via a methodology called \textit{topological data analysis}. The main tool within topological data analysis is \textit{persistent homology}; this is a means of measuring the shape of data, from the homology of a simplicial complex, calculated over a range of values.
The required background theory and a method of computing persistent homology is presented here, with applications specific to structural health monitoring. These results allow for topological inference and the ability to deduce features in higher-dimensional data, that might otherwise be overlooked. 

A simplicial complex is constructed for data for a given distance parameter. This complex encodes information about the local proximity of data points. A singular homology value can be calculated from this simplicial complex. Extending this idea, the distance parameter is given for a range of values, and the homology is calculated over this range. The persistent homology is a representation of how the homological features of the data persist over this interval. The result is characteristic to the data. A method that allows for the comparison of the persistent homology for different data sets is also discussed. 

\textbf{Key words: Topological data analysis; Persistent homology; Simplicial complex.}

\section{Introduction}
\label{sec:intro}
Topological methods are very rarely used in structural health monitoring (SHM), or indeed in structural dynamics generally, especially when considering the structure and topology of observed data. Topological methods can provide a way of proposing new metrics and methods of scrutinising data, the most rudimentary and most powerful of which, \textit{persistence homology}, will be discussed in this paper.

SHM has been dominated by the insurgence of machine learning since its introduction to the field \cite{farrar2007introduction,farrar2012structural}. \textit{Topological Data Analysis} (TDA) aims to work alongside the vast work already established to provide a new light into ways that data can be analysed. Previous machine learning operations in SHM have never distinctly considered the data \textit{topology}. This paper aims to bring to light the potential importance of TDA in SHM. TDA generalises well into higher dimensions, so an insight into higher-dimensional data structures can be obtained.

TDA has previously proven useful in other research areas, such as medicine, where virus evolution has been tracked. TDA has also been used to sequence viruses such as Influenza A and HIV, along with many other uses such as clustering diabetes types. In economics, TDA has been used to identify topological patterns in multi-dimensional time series data. It is thought that topological data analysis might be able to give early warning signs of imminent market crashes.

The contents of this paper aim to walk through the process of performing TDA. A brief outline of the methodology is given here, to construct a geometric object called a \textit{simplicial complex} from data. A simplicial complex can be thought of as a higher-dimensional analogue of a graph. The vertices of the simplicial complexes are the data points, and the connections between the vertices are less than a prescribed threshold distance. Therein, the simplicial complex encodes information about connection between the vertices. This geometric object has now attributed shape to the data. 
Following from the simplicial complex, this information can be manipulated in order to output algebraic groups that will capture information about the shape and structure of the connection in the simplicial complex, most specifically about the number of $k-$dimensional holes that will be found within the assumed manifold underlying the data; these groups are called the \textit{homology groups}. A generalisation of the homology will be discussed with respect to pioneering work by Edelsbrunner \cite{edelsbrunner2010computational, edelsbrunner2000topological} where the homology can be considered over a range of simplicial complexes; this is called the \textit{persistent homology}, aptly named as this uncovers how the homology persists over a range of values.

The layout of the paper is as follows: Section 2 give the necessary definitions, listed at the beginning to provide the required mathematics that is not usually common to engineering. Section 3 will be devoted to persistent homology, its significance, provide an intuitive understanding and show two common forms of how it can be displayed. Section 4 will show the topology of some known shapes and how the persistent homology links to these. Section 5 will introduce an interesting engineering specific example which will show how this theory can be used. Following this, the paper concludes.

\section{Background Theory and Definitions}
\label{sec:background}
At the expense of labouring multiple required definitions, the theory pertinent to the paper topic will quickly follow. For more details, the following works can be consulted \cite{rotman2012introduction, schutz1980geometrical, genomics, nash1988topology, ghrist2014EAT}.
\subsection{Algebra}
\begin{defn}
An \textit{equivalence relation}, denoted by $\sim$, is a binary operation that is \textit{reflexive} $(a \sim a)$, \textit{symmetric} (if $a \sim b$ then $b \sim a$) and \textit{transitive} (if $a \sim b$ and $b \sim c$ then $a \sim c$). The equivalence relation provides a partition of a set into elements that share a common property.
\end{defn}

\begin{defn}
\textit{Groups} are an extension of the concept of set, to include a binary operator. This set-operator pair is written $G=\{ S,\circ \}$. There are four necessary axioms associated with a group, these are: \textit{closure}, \textit{associativity}, \textit{existence of a unique identity}, and \textit{the existence of an inverse} for every element in the group.
\end{defn}

\begin{defn}
A \textit{subgroup}, $H$, of a group, $G$, is a subset such that $H \subseteq G$, and $H$ satisfies all of the group axioms. The operation and identity of the subgroup is inherited from the parent group.
\label{def: subgroup1}
\end{defn}

\begin{defn}
A group is called an \textit{abelian group}, if the result of the operation is independent of the ordering of elements. This is the case when the operation is \textit{commutative}. $x \circ y=y \circ x, \ \forall x,y \in G$.
\end{defn}

\begin{defn}
Let $G$ be a group and $H$ a subgroup of $G$. A \textit{left coset} of $H$ in $G$ is a subset of $G$ of the form $gH = \{gh\ |\ h \in H\}$ for some $g \in G$, the set of left cosets of $H$ in $G$ is written $G/H$. Similarly, a \textit{right coset} of $H$ in $G$ is a subset of $G$ of the form $Hg = \{hg\ |\ h \in H\}$ for some $g \in G$, and the set of right cosets is written $H\backslash G$ \cite{rotman2012introduction}. 
\end{defn}

\begin{defn}
A subgroup $H$ of $G$ is said to be \textit{normal} if $gH=Hg, \ \forall g \in G$, this is written $H \triangleleft G$. In this case, the set of cosets $G/H$ and $H\backslash G$ are the same by definition. The set of cosets for normal subgroups is called the \textit{quotient group}.
\end{defn}

\begin{defn}
A group \textit{homomorphism} is a \textit{mapping} between two groups, $h:G_1 \rightarrow G_2$. The mapping is considered a homomorphism if the identity element in $G_1$ is mapped to the identity element in $G_2$, and the group operation distributes over the homomorphism, $h(u*v)=h(u) \circ h(v),\forall u,v \in G_1$.
\end{defn}

\begin{defn}
Given a group homomorphism $h:G_1 \rightarrow G_2$, the \textit{kernel} of $h$, $\text{ker}(h) \subset G_1$, is the set of elements $x$ such that $h(x)=e$, where $e$ is the identity element. The \textit{image} of $h$, $\text{im}(h) \subset G_2$, is the set of elements $y$ such that $y=h(x)$ for some $x$.
\end{defn}

\begin{defn}
A \textit{metric space} is defined by a pair $(X,\partial_X)$. $X$ refers to the set where the elements of the metric space live and $\partial_X$ is the associated metric or distance function between two points in $X$. For a space to qualify as a metric space, the following criteria must be true: $\partial(x,y) \geq 0$, $\partial(x,y)=\partial(y,x)$, and $\partial(x,z) \leq \partial(x,y)+\partial(y,z)$.
\end{defn}

\begin{defn}
An \textit{open ball} is defined as $B_\epsilon(x) = \{y \in X \ | \ \partial_X(y,x)<\epsilon \}$. This encloses a space around the point $x$, where all points enclosed are less than the distance $\epsilon$ from the point $x$. This space is often referred to as the \textit{$\epsilon-$neighbourhood} of $x$.
\end{defn}

\begin{defn}
A \textit{topological space} is represented by a pair, $(X,\mathcal T)$, where $X$ is the set of all the elements and $\mathcal T$ is a collection of subsets, referred to as the \textit{topology}. The \textit{open sets} in $\mathcal T$ must satisfy the following axioms: The set of elements and the empty set are elements of the topology, any union of sets in $\mathcal T$ is also an element of $\mathcal T$, the intersection of a finite collection of elements of $\mathcal T$ is an element of $\mathcal T$.
\end{defn}

\begin{defn}
Given two topological spaces, $(X,\mathcal{T}_X)$ and $(Y,\mathcal{T}_Y)$. These spaces are said to be \textit{homeomorphic} if there exists a bijective continuous map between them.
\end{defn}

Now that a topological space has been discussed, when taking homeomorphisms between spaces there are quantities characteristic of that structure called \textit{topological invariants}; these do not change under homeomorphisms. It will be the main purpose of this review to develop a method of calculating a topological invariant called the \textit{persistent homology}.

\subsection{Manifolds}
Manifolds are continuous surfaces from which the data are assumed to be sampled. Data that are observed are assumed to lie on the surface of a manifold. By understanding the topology of the sampled data points, it is the aim of TDA to extract topological information about the underlying manifold from the sampled data. The manifold is unknown prior to analysis and persistent homology will identify features within the manifold over a range of length scales. Thereby, understanding the shape of the sampled data, it is the conjecture of TDA that the shape of the manifold is also understood. Formally, a manifold is a space that is locally homeomorphic to some $n-$dimensional Euclidean space, $\mathbb{R}^n$ \cite{schutz1980geometrical}.  
\subsection{Simplicial Complex}
Simplicial complexes are a way of representing data sampled from a manifold; they can be thought of as higher-dimensional analogues of graphs, giving a way of encoding connections between vertices. The dimension that a simplicial complex can capture is restricted by the number of points that are fully connected.   Simplicial complexes can be manipulated to output the homology for the data, and following this, the persistent homology.

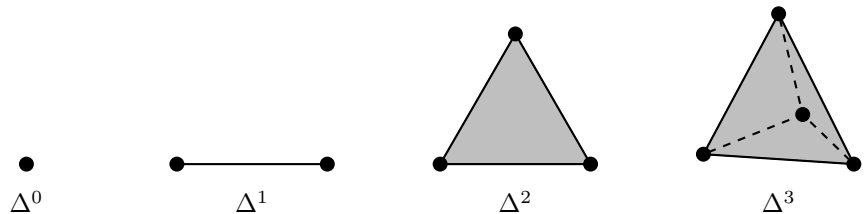
\begin{figure}[H]
    \centering
    
\begin{tikzpicture}[ele/.style={fill=black,circle,minimum width=3pt,inner sep=2pt}]
\node[ele] (a1) at (0,0) {};
\node at (0,-0.5) {$\Delta^0$};

\begin{scope}[shift={(2,0)}]
    \node[ele] (a1) at (0,0) {};    
    \node[ele] (a2) at (2,0) {};
    \draw[-,thick,shorten <=2pt,shorten >=2pt] (a1.center) -- (a2.center);
    \node at (1,-0.5) {$\Delta^1$};
\end{scope}

\begin{scope}[shift={(5.5,0)}]
    \node[ele] (a1) at (2,0) {};    
    \node[ele] (a2) at (0,0) {};
    \node[ele] (a3) at (1,1.732) {};

    \draw[-,thick,shorten <=2pt,shorten >=2pt] (a1.center) -- (a2.center);
    \draw[-,thick,shorten <=2pt,shorten >=2pt] (a1.center) -- (a3.center);
    \draw[-,thick,shorten <=2pt,shorten >=2pt] (a2.center) -- (a3.center);

    \begin{scope}[on background layer]
        \path [fill=lightgray,draw] (a1.center) to (a2.center) to (a3.center) to (a1.center);
    \end{scope}

    \node at (1,-0.5) {$\Delta^2$};
\end{scope}

\begin{scope}[shift={(9,0)}]
    \node[ele] (a1) at (2,0) {};    
    \node[ele] (a2) at (0,0.132) {};
    \node[ele] (a3) at (1.32,0.66) {};
    \node[ele] (a4) at (1,2) {};
    
    \draw[-,thick,shorten <=2pt,shorten >=2pt] (a1.center) -- (a2.center);
    \draw[-,thick,shorten <=2pt,shorten >=2pt] (a1.center) -- (a4.center);
    \draw[-,thick,shorten <=2pt,shorten >=2pt] (a2.center) -- (a4.center);

    \draw[-,thick,shorten <=2pt,shorten >=2pt, dashed] (a1.center) -- (a3.center);
    \draw[-,thick,shorten <=2pt,shorten >=2pt, dashed] (a2.center) -- (a3.center);
    \draw[-,thick,shorten <=2pt,shorten >=2pt, dashed] (a4.center) -- (a3.center);
    
    \begin{scope}[on background layer]
        \path [fill=lightgray,draw] (a1.center) to (a2.center) to (a4.center) to (a1.center);
    \end{scope} 
\node at (1,-0.5) {$\Delta^3$};
\end{scope}
\end{tikzpicture}
    \caption{The first four simplices.}
    \label{fig:simplices}
\end{figure}

A simplicial complex is a structure made up of fundamental building blocks called \textit{simplices}. The first four simplices can be seen in Figure \ref{fig:simplices}. Each vertex in the simplex is fully connected to all the other vertices and the space enclosed by the vertices is part of that simplex. 

There are many ways to construct a simplicial complex from point data. For simplification, only one method will be discussed within this paper, the \textit{Vietoris-Rips} complex. The Vietoris-Rips (VR) complex can be constructed for point data to output a corresponding complex according to the rules $VR_\epsilon(X,\partial_X)$: let $(X,\partial_X)$ be a finite metric space and $\epsilon > 0$ be a fixed value. The abstract simplicial complex is determined by the rules \cite{chambers2010vietoris}:
\begin{enumerate}
    \item The vertices, $v \in X$, form the vertices in $VR_\epsilon(X,\partial_X)$.
    \item A $k-$simplex is formed when $\partial_X(v_i,v_j) \leq 2\epsilon, \ \forall i,j \leq k$ for some $\varepsilon > 0$.
\end{enumerate}

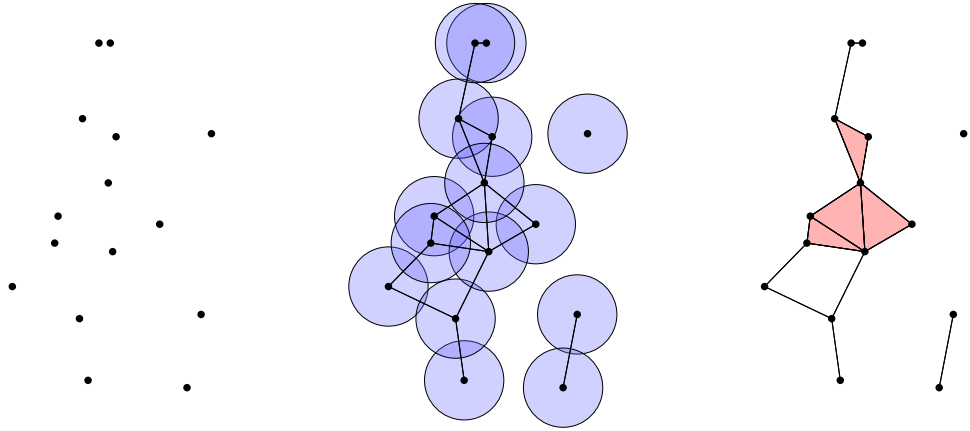
\begin{figure}[H]
    \centering
    \begin{tikzpicture}[ele/.style={fill=black,circle,minimum width=2pt,inner sep=1pt}, scale = 1]

    \node[ele] (a) at (2.8209087137275, 1.18093727264829) {};
    \node[ele] (b) at (2.95833231728538,	3.58594837328201) {};
    \node[ele] (c) at (0.311441915093690,	1.55228774215654) {};
    \node[ele] (d) at (1.69060866126193,	3.54516924716448) {};
    \node[ele] (e) at (1.31799287927655,	0.304557156170160) {};
    \node[ele] (f) at (1.61347687509477,	4.79163906483583) {};
    \node[ele] (g) at (2.27043233468693,	2.38135190156370) {};
    \node[ele] (h) at (1.20381702801476,	1.12542610389225) {};
    \node[ele] (i) at (0.919304590945695,	2.48862221210784) {};
    \node[ele] (j) at (1.64512903092519,	2.01677357394698) {};
    \node[ele] (k) at (1.24383442716840,	3.78580579686650) {};
    \node[ele] (l) at (1.58646607060451,	2.93089871716034) {};
    \node[ele] (m) at (0.874518089942222,	2.13017381291781) {};
    \node[ele] (n) at (1.46214937616206,	4.79124963787566) {};
    \node[ele] (o) at (2.63312304723148,	0.206799557528458) {};

\begin{scope}[elem/.style={circle, draw = black, line width = 0.2pt, fill = blue!60,fill opacity=0.3, minimum width = 30pt}, shift={(5,0)}]

\node[elem] (a) at (2.8209087137275, 1.18093727264829) {};
\node[elem] (b) at (2.95833231728538,   3.58594837328201) {};
\node[elem] (c) at (0.311441915093690,  1.55228774215654) {};
\node[elem] (d) at (1.69060866126193,   3.54516924716448) {};
\node[elem] (e) at (1.31799287927655,   0.304557156170160) {};
\node[elem] (f) at (1.61347687509477,   4.79163906483583) {};
\node[elem] (g) at (2.27043233468693,   2.38135190156370) {};
\node[elem] (h) at (1.20381702801476,   1.12542610389225) {};
\node[elem] (i) at (0.919304590945695,  2.48862221210784) {};
\node[elem] (j) at (1.64512903092519,   2.01677357394698) {};
\node[elem] (k) at (1.24383442716840,   3.78580579686650) {};
\node[elem] (l) at (1.58646607060451,   2.93089871716034) {};
\node[elem] (m) at (0.874518089942222,  2.13017381291781) {};
\node[elem] (n) at (1.46214937616206,   4.79124963787566) {};
\node[elem] (o) at (2.63312304723148,   0.206799557528458) {};

\node[ele] (a1) at (2.8209087137275, 1.18093727264829) {};
\node[ele] (b1) at (2.95833231728538,    3.58594837328201) {};
\node[ele] (c1) at (0.311441915093690,   1.55228774215654) {};
\node[ele] (d1) at (1.69060866126193,    3.54516924716448) {};
\node[ele] (e1) at (1.31799287927655,    0.304557156170160) {};
\node[ele] (f1) at (1.61347687509477,    4.79163906483583) {};
\node[ele] (g1) at (2.27043233468693,    2.38135190156370) {};
\node[ele] (h1) at (1.20381702801476,    1.12542610389225) {};
\node[ele] (i1) at (0.919304590945695,   2.48862221210784) {};
\node[ele] (j1) at (1.64512903092519,    2.01677357394698) {};
\node[ele] (k1) at (1.24383442716840,    3.78580579686650) {};
\node[ele] (l1) at (1.58646607060451,    2.93089871716034) {};
\node[ele] (m1) at (0.874518089942222,   2.13017381291781) {};
\node[ele] (n1) at (1.46214937616206,    4.79124963787566) {};
\node[ele] (o1) at (2.63312304723148,    0.206799557528458) {};
\foreach \firstnode in {a, b, c, d, e, f, g, h, i, j, k, l, m, n, o}{%
   \foreach \secondnode in {a, b, c, d, e, f, g, h, i, j, k, l, m, n, o}{%
   \setveclength{\mydist}{\firstnode}{\secondnode}
   \pgfmathparse{\mydist < 30 ? int(1) : int(0)}
   \ifnum\pgfmathresult=1
     \draw (\firstnode.center) -- (\secondnode.center);
   \fi
  }
}
\end{scope}

\begin{scope}[shift={(10,0)}]
    \node[ele] (a) at (2.8209087137275, 1.18093727264829) {};
    \node[ele] (b) at (2.95833231728538,	3.58594837328201) {};
    \node[ele] (c) at (0.311441915093690,	1.55228774215654) {};
    \node[ele] (d) at (1.69060866126193,	3.54516924716448) {};
    \node[ele] (e) at (1.31799287927655,	0.304557156170160) {};
    \node[ele] (f) at (1.61347687509477,	4.79163906483583) {};
    \node[ele] (g) at (2.27043233468693,	2.38135190156370) {};
    \node[ele] (h) at (1.20381702801476,	1.12542610389225) {};
    \node[ele] (i) at (0.919304590945695,	2.48862221210784) {};
    \node[ele] (j) at (1.64512903092519,	2.01677357394698) {};
    \node[ele] (k) at (1.24383442716840,	3.78580579686650) {};
    \node[ele] (l) at (1.58646607060451,	2.93089871716034) {};
    \node[ele] (m) at (0.874518089942222,	2.13017381291781) {};
    \node[ele] (n) at (1.46214937616206,	4.79124963787566) {};
    \node[ele] (o) at (2.63312304723148,	0.206799557528458) {};
    \foreach \firstnode in {a, b, c, d, e, f, g, h, i, j, k, l, m, n, o}{%
   \foreach \secondnode in {a, b, c, d, e, f, g, h, i, j, k, l, m, n, o}{%
   \setveclength{\mydist}{\firstnode}{\secondnode}
   \pgfmathparse{\mydist < 30 ? int(1) : int(0)}
   \ifnum\pgfmathresult=1
     \draw (\firstnode.center) -- (\secondnode.center);
   \fi
  }
}
\begin{pgfonlayer}{background}
\path[fill=red!30, draw] (g.center) to (j.center) to (l.center) to (g.center); 
\path[fill=red!30, draw] (d.center) to (k.center) to (l.center) to (d.center); 
\path[fill=red!30, draw] (m.center) to (i.center) to (j.center) to (m.center); 
\path[fill=red!30, draw] (l.center) to (j.center) to (i.center) to (l.center); 
\end{pgfonlayer}

\end{scope}

\end{tikzpicture}
    \caption{The process of constructing a VR complex from some data.}
    \label{fig: Rips}
\end{figure}

This process of constructing a simplicial complex for some data attributes shape to the data, by defining connections between arbitrarily close points. Before moving on to homology, which is an algebraic construct, a method of converting the geometric simplicial complex into an algebraic simplicial complex must first be considered. An abstract simplicial complex represents a method of defining the geometric connections, so that algebra can then be applied to the shape.
\begin{defn}
An \textit{abstract simplicial complex}, $K$,  consists of a set of vertices, $\text{vert}(K)$ and a set of abstract simplices, $\text{simp}(K)$, such that \cite{chazal2017introduction}:
\begin{enumerate}
    \item Every simplex $\Delta^k \in \text{simp}(K)$ is a non-empty subset of $\text{vert}(K)$, or the simplices are a union of the vertices.
    \item For all vertices, $v \in \text{vert}(K)$, there is also an abstract simplex $\{ v \} \in \text{simp}(K)$.
    \item For every non-empty abstract simplex, $\Delta^k$, a non-empty subset of $\Delta^k$ is also an abstract simplex. This is referred to as a \textit{face} of $\Delta^k$.
    \item For an abstract simplex $\Delta^k \in \text{vert}(K)$, the dimension of $\Delta^k$ is $\text{dim}(\Delta^k) = |\Delta^k|-1$, where $|.|$ denotes the number of elements in the set.
\end{enumerate}
\label{def:simpcomp}
\end{defn}
This notion of an abstract simplicial complex is pivotal to TDA. Abstract simplicial complexes allow for the mapping between categories: a mapping from a geometric realisation of a simplicial complex which is embedded in Euclidean space, to an abstract simplicial complex. The abstract simplicial complex form can be algorithmically analysed with computer packages to output results for that simplicial complex.

\subsection{Homology}
The homology groups, $H_k(X)$ are invariants for the data set, $X$, where $k$ refers to the dimension of the feature of the homology group. The homology groups encode information about the number of $k-$dimensional holes in the data.

\begin{defn}
The \textit{orientation} of the vertices of a simplex, $\Delta^k$, is an equivalence class of orderings of the vertices under the equivalence relation that two orderings are the same if they differ by an \textit{even permutation}. An even permutation is one that can be expressed as a composition of even permutations \cite{genomics}. There are only two possible orientations, a positive and a negative.
\end{defn}

\begin{defn} \label{def: chaingroup}
To each standard simplex $\Delta^k_i$ of a simplicial complex $K$, an abelian group, $C_k(K)$, called a \textit{chain group} is associated to it. The $k^{th}$ chain group is the set of all $k-$simplices in the simplicial complex $K$ \cite{ghrist2018homological}.

Now, consider a $K$ containing $l_k$ of any standard simplex $\Delta^k$ for all values of $k$. The \textit{$k-$chain group} of $K$, $C_k(K)$ is the free abelian group generated by the oriented $k-$simplices of $K$. This means any element $c_k \in C_k(K)$ can be thought of in abstract as, $$c_k = \sum_{i=1}^{l_k}f_i\Delta_i^k, \ f_i \in \mathbb{Z}$$
where the following criteria are satisfied \cite{nash1988topology}:
\begin{enumerate}
    \item A negation of simplices. $\Delta_i^k+ (-\Delta_i^k) = 0, \ \forall i,k$.
    \item A linearity over the elements. $\displaystyle \sum_{i=1}^{l_k}f_i\Delta_i^k + \sum_{i=1}^{l_k}g_i\Delta_i^k = \sum_{i=1}^{l_k}(f_i+g_i)\Delta_i^k$ where $f_i,g_i\in \mathbb{Z}$.
\end{enumerate}
\end{defn}

Following the previous two definitions, a boundary operator can be formulated. The boundary operator provides an algebraic formulation for the exterior region of a simplicial complex. Now that simplicial complexes are being represented by chain groups, the problem has become entirely algebraic.

\begin{defn}\label{def: boundaryvertremoved}
The \textit{boundary operator}, $\partial_k$, maps between chain groups, 
$$\partial_k \colon C_k(K)\rightarrow C_{k-1}(K)$$
Given an oriented simplex $\Delta^k = [v_0,...,v_k]$, a positive sign is assigned to every member of the even permutation class of $\Delta^k$ and a negative sign to every member of the odd permutation class. The boundary operator must now obey the rules \cite{nash1988topology}:
\begin{enumerate}
    \item For an oriented simplex, 
    $$\partial \Delta^k = \sum_{i=0}^k(-1)^i[v_0,...,\hat{v}_i,...,v_k] $$
    where $[v_0,...,\hat{v}_i,...,v_k]$ represents the face of the simplex with the $i^{th}$ vertex omitted. Note that every successive omission changes the orientation of the face.
    \item Thinking of a simplicial complex $K$, as the abstract sum of all the standard simplices required to construct it, $K = \sum_{i,k}\Delta^k_i, \ \Delta^k_i \subset K$. The boundary operator is a linear function over all the simplices in the complex.$$\partial(K) = \partial\big(\sum_{i,k}\Delta^k_i \big) = \sum_{i,k}\partial (\Delta^k_i) , \ \Delta^k_i \subseteq K$$
\end{enumerate}
\end{defn}

\begin{defn}
A mapping between successive chain groups with the boundary operator,
$$... \rightarrow C_k \xrightarrow{\partial_{k}} C_{k-1} \xrightarrow{\partial_{k-1}} ... \xrightarrow{\partial_{2}} C_1 \xrightarrow{\partial_{1}} C_0 \xrightarrow{\partial_{0}} 0$$
is called a \textit{chain complex} \cite{ghrist2018homological}.
\end{defn}
For the purposes of calculations, it is said that $C_i(K)=0$ for $i>\text{dim}(K)$. It is also the case that $C_i(K)=0$ for $i<0$. In both these cases, 0 represents the zero group.
\begin{defn}\label{def: kcycle}
For a simplicial complex $K$, elements of the chain group, $z_k\in C_k(K)$, are called \textit{$k-$cycles} if $\partial z_k = 0$. The group of $k-$cycles, $Z_k(K)$, is given by the kernel of the boundary map,
$$Z_k(K) = \text{ker}(\partial_k \colon C_k(K) \rightarrow C_{k-1}(K))= \{z_k \in C_k \colon \partial_k z_k = 0\}$$
and $Z_k(K)$ is a subgroup of $C_k(K)$ \cite{boissonnat2018geometric}.
\end{defn}
\begin{defn}\label{def: kboundary}
For a simplicial complex $K$, elements of the chain group, $b_k\in C_k(K)$, are called \textit{$k-$boundaries} if there exists a $(k+1)-$chain group, $C_{k+1}(K)$, such that $\partial C_{k+1}(K)=b_k$. The group of $k-$boundaries, $B_k(K)$, is given by,
$$B_k(K) = \text{im}(\partial_{k+1} \colon C_{k+1}(K) \rightarrow C_{k}(K))= \{b_k \in C_k \colon \exists b_k' \in C_{k+1}, b_k = \partial b_k'\}$$
and $B_k(K)$ is a subgroup of $C_k(K)$ \cite{boissonnat2018geometric}.
\label{def:kbound}
\end{defn}

The interesting part of this analysis, which fundamentally results in the homology groups, is that when the composition of two boundary operations is analysed,  applying it twice gives a result of 0. That is the boundary of a boundary is empty.
\begin{lemm}
$$\partial_{k-1} \circ \partial_k = 0$$
\label{lemm: boundaryzero}
This is a very important result, a proof is supplied in many texts \cite{schutz1980geometrical,nash1988topology} 
\end{lemm}

\begin{corr}
A very important result follows on from Lemma \ref{lemm: boundaryzero}. For a simplicial complex, $K$, any element of the boundary group $b_k \in B_k(K)$ has the property $\partial_k b_k = 0$. Therefore, $B_k(K) \subseteq Z_k(K)$ where $Z_k(K)$ is the group of $k-$cycles. Since both $Z_k(K)$ and $B_k(K)$ are abelian, a property inherited by being subgroups of $C_k(K)$. $B_k(K)$ is a normal subgroup of $Z_k(K)$, and therefore it can divide the parent group.
\label{corr: homologydivides}
\end{corr}
\begin{figure}[H]
    \centering
    \begin{tikzpicture}
\draw[fill=cyan!10] (0,1) ellipse (1.5cm and 2cm);
\draw[fill=cyan!30] (0,0.5) ellipse (1cm and 1.5cm);
\draw[fill=cyan!50] (0,0) ellipse (0.5cm and 1cm);
\node at (0,0.25) {$B_{p+1}$};
\node at (0,1.5) {$Z_{p+1}$};
\node at (0,2.5) {$C_{p+1}$};

\draw[fill=cyan!10] (5,1) ellipse (1.5cm and 2cm);
\draw[fill=cyan!30] (5,0.5) ellipse (1cm and 1.5cm);
\draw[fill=cyan!50] (5,0) ellipse (0.5cm and 1cm);
\node at (5,0.25) {$B_{p}$};
\node at (5,1.5) {$Z_{p}$};
\node at (5,2.5) {$C_{p}$};

\draw[fill=cyan!10] (10,1) ellipse (1.5cm and 2cm);
\draw[fill=cyan!30] (10,0.5) ellipse (1cm and 1.5cm);
\draw[fill=cyan!50] (10,0) ellipse (0.5cm and 1cm);
\node at (10,0.25) {$B_{p-1}$};
\node at (10,1.5) {$Z_{p-1}$};
\node at (10,2.5) {$C_{p-1}$};

\draw (-2.5,-1) [dashed]-- (0,-1);
\draw (0,-1) -- (10,-1);
\draw (10,-1) [dashed]-- (12.5,-1);

\draw (0.5,1.8) -- (5,-1);
\draw (0.573,2.85) -- (5.07,0.99);
\draw (5.5,1.8) -- (10,-1);
\draw (5.573,2.85) -- (10.07,0.99);

\draw[->] (-3,-0.5) to node[midway,above] {$\partial_{p+2}$} (-2,-0.5);
\draw[->] (2,-0.5) to node[midway,above] {$\partial_{p+1}$} (3,-0.5);
\draw[->] (7,-0.5) to node[midway,above] {$\partial_{p}$} (8,-0.5);
\draw[->] (12,-0.5) to node[midway,above] {$\partial_{p-1}$} (13,-0.5);

\draw[fill = white]  (0,-1) circle (0.1);
\draw[fill = white] (5,-1) circle (0.1);
\draw[fill = white] (10,-1) circle (0.1);

\node[below] at (0,-1.1) {$0$};
\node[below] at (5,-1.1) {$0$};
\node[below] at (10,-1.1) {$0$};
\end{tikzpicture}
    \caption{A chain complex with boundary maps between the chain groups.}    \label{fig:chaincomp}
\end{figure}
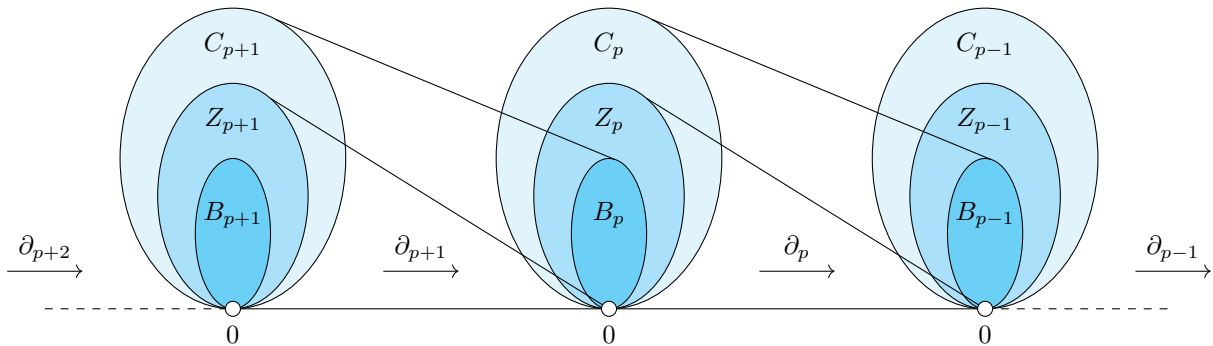
Corollary \ref{corr: homologydivides} brings together a lot of the background theory that has been covered up until this point. It is beneficial to break this down into bite-size chunks. The first stage is understanding that $\partial_kb_k=0$, as this is essentially the application of two boundary operators, as the element $b_k$ is formed from the boundary operation of a previous chain element. This can be rewritten as $(\partial_k \circ \partial_{k+1})c_{k+1}$, then using Lemma \ref{lemm: boundaryzero} it follows that this result is equal to zero. If all the elements in $B_k(K)$ are mapped to zero after the application of $\partial_k$ then it must be a subset of the group that is mapped to zero, by Definition \ref{def: kcycle} this is $Z_k(K)$. From Definitions \ref{def: kcycle} and \ref{def: kboundary}, it is known that $Z_k(K) \subset C_k(K)$ and $B_k(K) \subset C_k(K)$. Since the chain groups $C_k(K)$ are defined to be abelian, then so are their subgroups, $B_k(K)$ and $Z_k(K)$ by Definition \ref{def: subgroup1}. It has also been proven that $B_k(K) \subset Z_k(K)$, and since $B_k(K)$ is abelian it is known that $B_k(K) \triangleleft Z_k(K)$ and therefore the set of cosets $Z_k(K)/B_k(K)$ is the quotient group.

\begin{defn}
The \textit{homology groups}, $H_k(K)$, are the quotient groups,
$$H_k = Z_k(K) / B_k(K)$$
\end{defn}
Setting aside all this mathematical rigour, the $k^{th}$ homology groups can simply be thought of as the cycles in $C_k$ that are not boundaries of the elements within $C_{k+1}$. That an element of $C_k$ is a cycle, means it encloses a $k-$dimensional region. The fact this is not a boundary means the interior is not part of the underlying space. This is where the idea of counting the $k-$dimensional holes springs from. A generalisation of the rule is that $H_k(K)$ counts the $k-$dimensional holes in the simplicial complex, $K$. $H_0(K)$ is the only real exception to this rule, as this encodes information about the number of path-connected components in $K$. $H_1(K)$ encodes information about $1D$ holes; these can be visualised as circular holes. $H_2(K)$ encodes information about $2D$ holes; these can be visualised as cavities. $H_k(K)$ encodes information about $kD$ holes.

\begin{defn}
$H_k(K)$ is a vector space and the elements are the \textit{homology classes} of $K$. The homology class of a cycle $z_k \in Z_k(K)$ is the coset $c_k + B_k(K) =\{ c_k + b_k \colon b_k \in B_k(K) \}$. Cycles are said to be \textit{homologous} if they are in the same homology class \cite{boissonnat2018geometric}. 
\end{defn}

\begin{defn}
The \textit{Betti number}, $\beta_k$, is the dimension of the $k^{th}$ homology group of a simplicial complex,
$$\beta_k = \text{dim}(H_k(K))$$
The Betti number represents the number of $k-$dimensional holes in a simplicial complex \cite{ghrist2014EAT}.
\end{defn}
Betti numbers will be the primary topological invariants used throughout this paper. Betti numbers will be used in persistent homology and they are vital in visualising spaces. Examples of Betti numbers and how they can be visualised are provided in Section \ref{sec:understandTDA}. Following this, Betti numbers will be used in the engineering examples in Section \ref{sec:application}.

\section{Persistent Homology}
\label{sec:persisthom}

\subsection{Understanding}
When the distance, $\epsilon$, is smaller than some feature scale, the properties of that feature can be captured. This means topological invariants that are described at a length scale greater than $\epsilon$ can be captured. A problem arises here, as usually the feature scale is not known prior to analysis. 

Obtaining the homology for a single value of $\epsilon$ provides very limited information, this notion is almost redundant, due to potential varying feature length scales in the manifold. For this reason, it is vital to consider what homological features persist as $\epsilon$ is varied. The goal of \textit{persistent homology} is to track the homology classes as $\epsilon$ is varied. This process of varying $\epsilon$ does not bias any disk size, as all are being considered. This process will give an initial value, $\epsilon_{\text{min}}$, where a specific homological feature comes to life and $\epsilon_{\text{max}}$, where the feature is no longer considered for that simplicial complex. This range of values $[\epsilon_{\text{min}},\epsilon_{\text{max}}]$ is called the \textit{persistence interval} for that homological feature. Each persistence interval is attributed a Betti number. Following this, the set of all persistence intervals is descriptive for that manifold, giving information about in which dimension a hole exists in the data and over what range of values it persists for. Regardless of triangulation of a simplicial complex, or construction method, the information obtained from persistence is the same; i.e. the persistence of a space is considered a topological invariant.

When varying $\epsilon$, the simplicial complex goes from $k$ disconnected vertices to a fully connected $(k-1)-$dimensional simplex (the filtration is likely to not be taken this far, as computational time would be unnecessarily large), where there are $k$ points sampled from a manifold. The features that persist over the longest relative change are more likely to represent features of the space. This procedure will then give a strong indication of the likelihood of the topological invariants associated to that space.

When calculating persistence, a minimum persistence interval length is specified. Given a persistence interval $[\epsilon_i,\epsilon_j]$, the length of the interval can be calculated by $|\epsilon_j-\epsilon_i|>l$ where $l$ represents the interval length threshold. Interval lengths shorter than $l$ will not be considered in the analysis;  this is useful as it could be argued that for very small persistence intervals, whether they even 'persist' or not.

The persistence intervals obtained can can be represented in two ways: barcodes or persistence diagrams, both having their merits.

\subsection{Diagrams}
\begin{figure}[H]
    \centering
    \includegraphics[width = 0.9\textwidth]{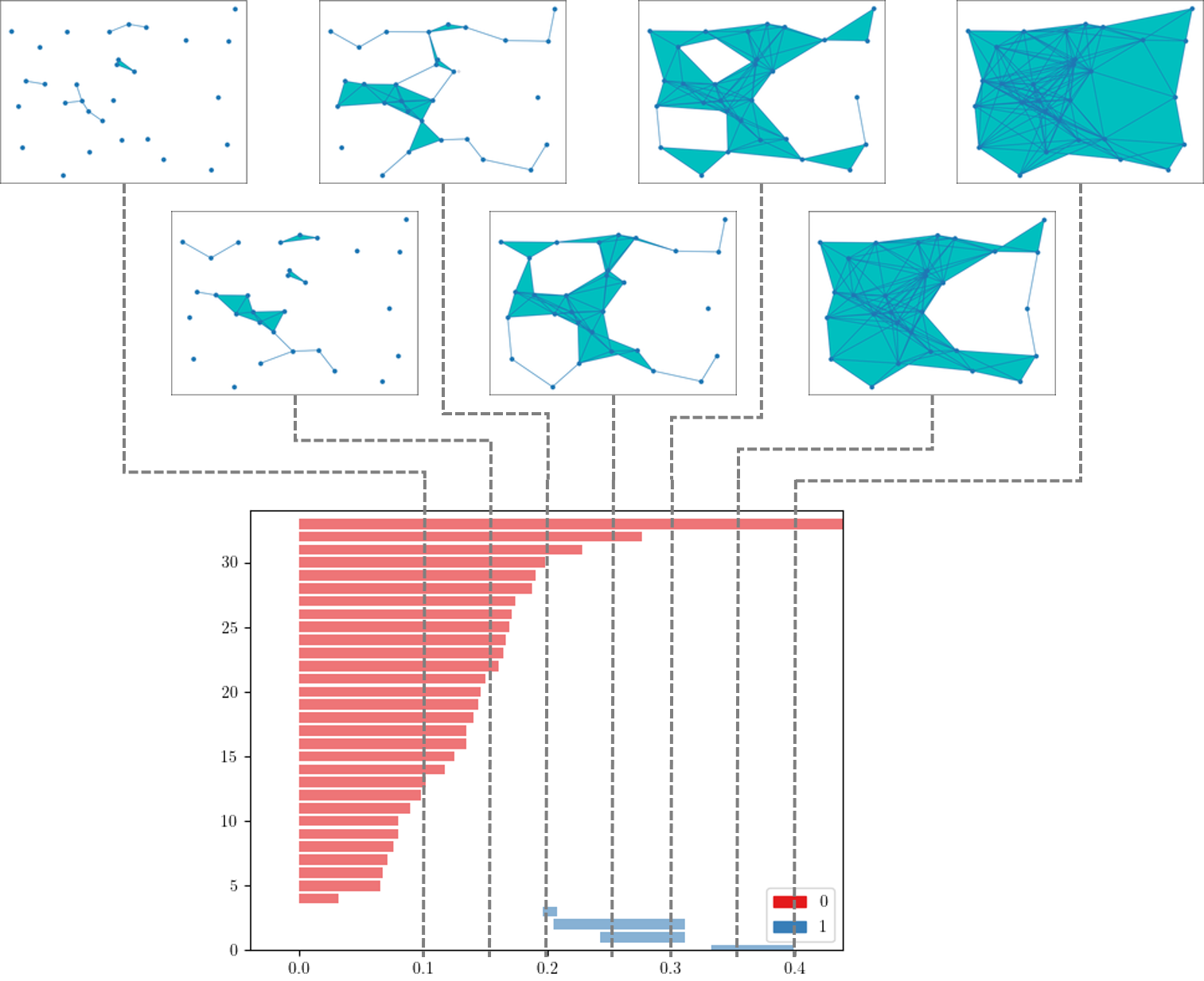}
    \caption{Persistence barcode with realisations showing which simplicial complexes are present at some values of $\varepsilon$.}
    \label{fig:persistencerealisations}
\end{figure}

In a barcode representation, the $x-$axis refers to the value of $\epsilon$. As $\epsilon$ increases, the barcode shows which features persist. The set of intervals are plotted with each interval beginning at $\epsilon_{\text{min}}$ and ending at $\epsilon_{\text{max}}$. The colour of the interval on the barcode refers to  the Betti number, $\beta_k$ \cite{ghrist2008barcodes}. The value of the $y-$axis can simply thought as an indexing of the intervals in the barcode. An example of a barcode can be seen in Figure \ref{fig:persistencerealisations}, with vertical dotted lines showing the intersections with the intervals, showing which features are present and the corresponding simplicial complex is found at the end of the dotted line. A few notable tricks to reading the barcodes are: the length of the interval represents how long the feature persists for. The length of the interval can be thought of as having a higher probability that this feature is characteristic of the manifold, as the longer the interval, the more prominent that feature is in the data.

\begin{figure}[H]
    \centering
    \includegraphics[width = 0.6\textwidth]{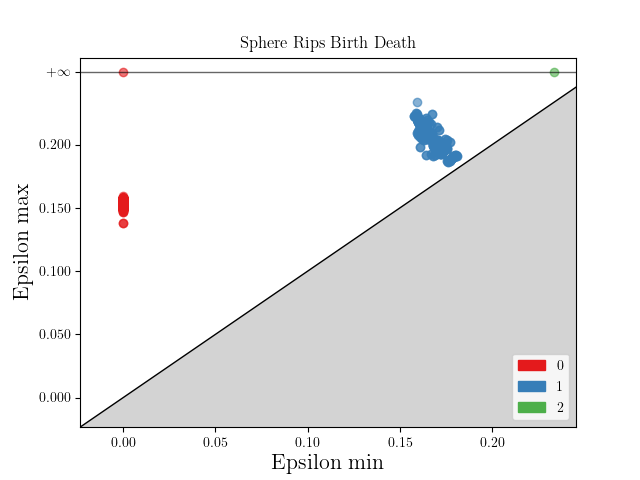}
    \caption{Birth-death diagram.}
    \label{fig:BDdia}
\end{figure}

The other method of visualising the set of persistence intervals is the persistence diagram, or \textit{birth-death} diagram. In this representation, $\epsilon_{\text{min}}$ is plotted on the $x-$axis and $\epsilon_{\text{max}}$ is plotted on the $y-$axis, with each interval represented by the point $(\epsilon_{\text{min}},\epsilon_{\text{min}})$. Intuitively, there is a line defined by $y=x$ below which the points will not be plotted; this line has the interpretation that the feature must first exist before it can die. Reading these diagrams is as intuitive as reading the barcodes; the vertical height of the point from the line $y=x$ is analogous to the length of the interval, that is the further a point is from the line $y=x$ the more the feature persists. An example of a birth-death diagram can be seen in Figure \ref{fig:BDdia}. The Birth-Death diagram shows a more structured method of comparing persistence intervals as the $y-$axis is not so arbitrary when compared to the barcodes.

On both the barcode and birth-death diagram, it can be seen that features persist to $\infty$. This has two meanings. The first being that there will always be a fully-connected simplex that persists to infinity. There will be a value of $\epsilon_{\text{fc}}$ that results in a fully-connected simplex where every vertex is connected to every other vertex. For values $\epsilon>\epsilon_{\text{fc}}$ the simplex will remain fully connected, and therefore this will continue to infinity. The second occurrence of this is when a feature persists past the value of $\epsilon$ used in the construction of the simplicial complex, as only smaller simplicial complexes of up to that value of $\epsilon$ are considered.

Each representation has their merits. For barcodes, it is easy to see which features a simplicial complex will have by drawing a vertical line at $x=\epsilon$, as can be seen in Figure \ref{fig:persistencerealisations}. The barcode representation shows repeated intervals, as these are new entries. This is not the case with the birth-death diagram, these intervals will displayed as the same point and overlap, therefore, information is lost in the birth-death diagram in this situation. Despite this, birth-death diagrams are less arbitrary and are displayed more compactly.

The space of barcodes forms a metric space; the distance between the barcodes is a measure of similarity of two barcodes. 

As manifolds can be represented by their barcodes, this notion of a metric space allows one to compare the similarity of manifolds. Metrics between barcodes are well established and the one used in this report is the \textit{$p-$Wasserstein distance}.

\begin{defn}
Given two barcodes $B_1$ and $B_2$. For $p>0$, the \textit{$p-$Wasserstein distance} is given by,
$$\partial_{W_p}(B_1,B_2) = \left(\text{inf} \sum_{Z\in B_1}d_\infty(Z, \phi(Z))^p \right)^{\frac{1}{p}}$$
where $\phi$ is a matching between $B_1$ and $B_2$ and $Z$ is a persistence interval in $B_1$ \cite{genomics}.
\end{defn}

\section{Understanding TDA}
\label{sec:understandTDA}
To help elaborate the point and highlight some features of TDA, two examples will be given in this section, where the points will be sampled from the manifold $S^2$ (the two-dimensional sphere embedded in $\mathbb{R}^3$) in two different ways. The first method used is by taking concentric circles in the plane perpendicular to the $z-$axis. The second method was to construct a \textit{Fibonacci spiral} around the sphere. The formula used for the first method to generate the point cloud (the set of samples from the manifold) for the embedding of the sphere was, 

\begin{minipage}{0.5\textwidth}
$$f(u,v) = \begin{cases}
    x=\cos u \sin v \\
    y=\sin u \cos v \\
    z=\cos v  \\
\end{cases}
 u \in [0,2\pi], \ v \in [0,\pi]$$
\end{minipage}
\hfill
\begin{minipage}{0.45\textwidth}
\includegraphics[width = 0.9\textwidth]{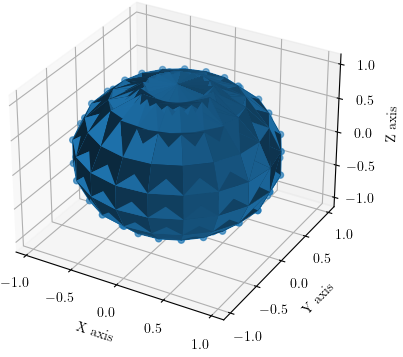}
\captionof{figure}{Sphere simplicial complex}
\null
\par\xdef\tpd{\the\prevdepth}
\end{minipage}

When constructing a VR complex, simplices are formed at all points within an open ball $B_\epsilon$; this means that connections will be formed to points past other points in the same direction. In this example, this is problematic as it results in the shape being over-connected at the poles if one wants the equatorial region to be filled with 2-simplices. This negatively impacts the computation time.
\begin{figure}[H]
    \centering
    \includegraphics[width = \textwidth]{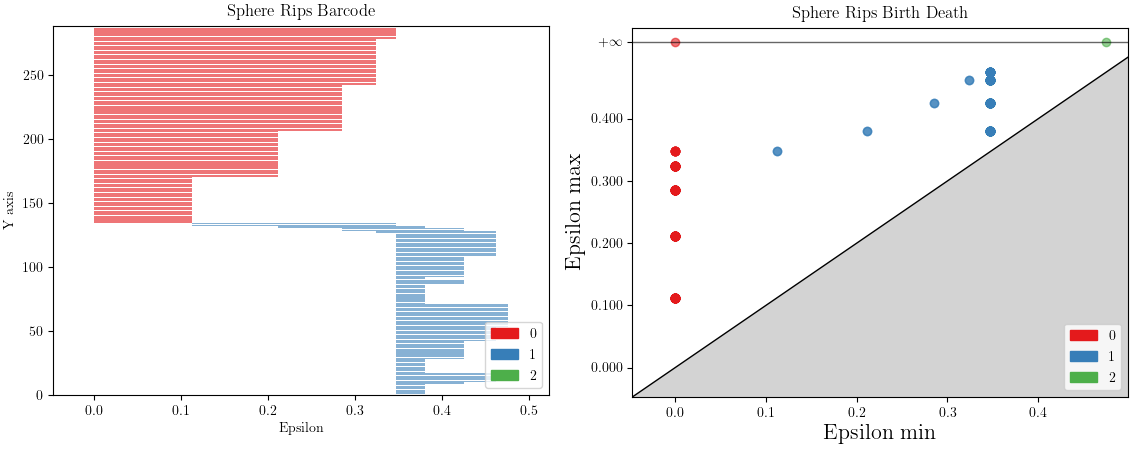}
    \caption{Sphere persistence representations.}
    \label{fig:Sph}
\end{figure}
The data and values for the given value of $\varepsilon$, for this sphere and the Fibonacci sphere is shown in Table \ref{tab:sphere}. It should be noted that the values of the Betti numbers in the table are specific to that value of epsilon. The diagrams contain a greater depth of information as the value of epsilon can be seen to be varied.

The Betti numbers in Table \ref{tab:sphere} are representative of the simplicial complex shown in Figure \ref{fig:Sph}. These are $\beta_0 = 1$, which says that the simplicial complex is fully connected with no disjoint parts. $\beta_1 = 0$, tells us that there are no holes in the simplicial complex. $\beta_2 = 1$ says that the simplicial complex enclose a volume, which is the inside of the sphere. It is important to note that the inside of the sphere is not part of the simplicial complex, but the space in which the simplicial complex is embedded.

The second case of sampling from the sphere, was done with \textit{Fibonacci spirals} around the sphere; this gives a much more uniform distribution of points. The formula used to generate the point cloud for the embedding of the sphere was,

\begin{minipage}{0.5\textwidth}
$$f(\phi,\theta) = \begin{cases}
    x=\cos \theta \sin \phi \\
    y=\sin \theta \sin \phi \\
    z=\cos \phi  \\
\end{cases}
\phi \in [0,\pi], \ \theta \in [0,n_p\pi(1+\sqrt{5})]$$
\end{minipage}
\hfill
\begin{minipage}{0.45\textwidth}
\includegraphics[width = 0.9\textwidth]{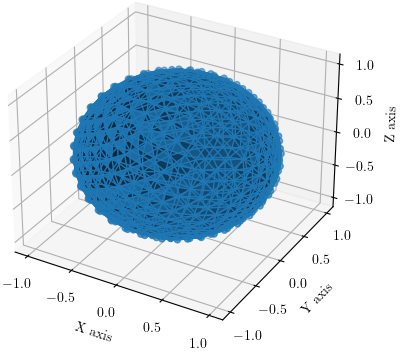}
\captionof{figure}{Fibonacci sphere simplicial complex.}

\null
\par\xdef\tpd{\the\prevdepth}
\end{minipage}

where $n_p$ is the number of points being sampled, in this case $n_p=500$.
\begin{figure}[H]
    \centering
    \includegraphics[width = \textwidth]{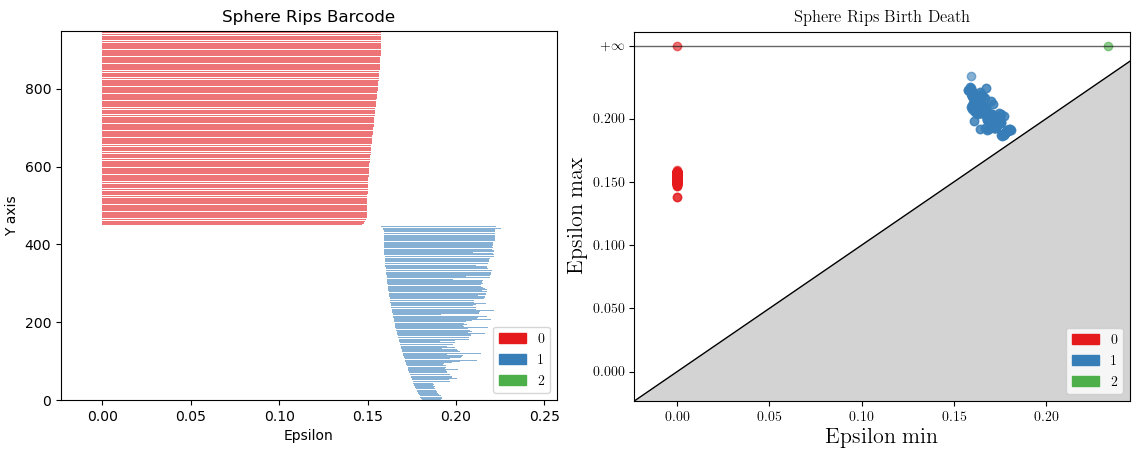}
    \caption{Fibonacci sphere persistence representations.}
    \label{fig:SphFib}
\end{figure}

From Table \ref{tab:sphere} it can be seen that the first scenario has fewer vertices, but still results in more simplices after constructing the simplicial complex. This is due to the poor distribution of points and the larger epsilon. The example goes to show that data cleansing steps and the correct model should be used to output meaningful and efficient results from TDA.

\begin{table}[H]
    \centering
    \begin{tabular}{|c|c|c|}
        \hline \rowcolor[HTML]{c0c0c0}
        Description   & Sphere & Fibonacci sphere  \\ \hline
        \rowcolor[HTML]{EFEFEF} 
        $\epsilon$      & $0.5$            & $0.25$            \\ \hline
        Number of vertices  & 200            & 500             \\ \hline
        \rowcolor[HTML]{EFEFEF} 
        Number of simplices  & 112094            & 4202             \\ \hline
        Betti Numbers & $[1,0,1]$           & $[1,0,1]$             \\ \hline
        \rowcolor[HTML]{EFEFEF} 
        Dimension     & 3               & 3                 \\ \hline
    \end{tabular}
        \caption{Topological properties for each sphere.}
        \label{tab:sphere}
\end{table}

When using the Wasserstein metric between the two data samples, a relatively small value of $\partial_{W_p}  = 2.673$ is obtained; this shows that the manifolds are relatively similar. This result is to be expected as the two spheres are sampled from the same manifold. When qualitatively analysing the persistence data, this backs up the case as the barcodes and birth-death diagrams, whilst not the same, do exhibit very similar features.

\section{Application}
\label{sec:application}
In this case, a four-dimensional surface plot will be considered. This example is defined by linear algebraic equations, where the variables form the axes of the higher-dimensional plots.

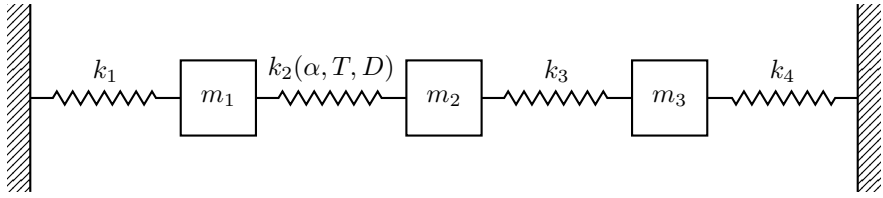
\begin{figure}[H]
    \centering
    \begin{tikzpicture}[every node/.style={outer sep=0pt},thick,
 mass/.style={draw,thick},
 spring/.style={thick,decorate,decoration={zigzag,pre length=0.3cm,post
 length=0.3cm,segment length=6}},
 ground/.style={fill,pattern=north east lines,draw=none,minimum
 width=0.75cm,minimum height=0.3cm},
 dampic/.pic={\fill[white] (-0.1,-0.3) rectangle (0.3,0.3);
 \draw (-0.3,0.3) -| (0.3,-0.3) -- (-0.3,-0.3);
 \draw[line width=1mm] (-0.1,-0.3) -- (-0.1,0.3);}]

  \node[mass,minimum width=1cm,minimum height=1cm] (m1) {$m_1$};
  \node[mass,minimum width=1cm,minimum height=1cm,right=2cm of
  m1] (m2) {$m_2$};
  \node[mass,minimum width=1cm,minimum height=1cm,right=2cm of
  m2] (m3) {$m_3$};
  \node[left=2cm of m1,ground,minimum width=3mm,minimum height=2.5cm] (g1){};
  \draw (g1.north east) -- (g1.south east);
  \node[right=2cm of m3,ground,minimum width=3mm,minimum height=2.5cm] (g2){};
  \draw (g2.north west) -- (g2.south west);

  \draw[spring] (g1.east) coordinate(aux)-- (m1.west|-aux) node[midway,above=1mm]{$k_1$};
  \draw[spring]  (m1.east|-aux) -- (m2.west|-aux) node[midway,above=1mm]{$k_2(\alpha,T,D)$};
  \draw[spring]  (m2.east|-aux) -- (m3.west|-aux) node[midway,above=1mm]{$k_3$};
  \draw[spring] (g2.west) coordinate(aux)-- (m3.east|-aux) node[midway,above=1mm]{$k_4$};

\end{tikzpicture}
    \caption{General mass-spring-damper system with variable $k_2$.}
    \label{fig:generalfbd}
\end{figure}

The system under analysis is a simple 3DOF system, with springs between masses and ground. A representation of a general case can be seen in Figure \ref{fig:generalfbd}. In this modelling scenario, it is assumed that the stiffness of the second spring, $k_2$ is a function of the coefficient of thermal expansion, $\alpha_2$, the temperature of the spring, $T_2$, and the damage present in the spring, $D_2$. It is the aim of this analysis to understand the topology of the manifold constructed when the natural frequency $\omega_i$ is calculated as a function of the variable parameters $\alpha_2, T_2, D_2$ and all other values are constant.
A simple relationship is assumed, where an increase in the temperature of a spring causes a reduction of the stiffness, by the relation $k_2(1-\alpha_2 T_2)$. It is also assumed that the damage coefficient also reduces the stiffness of the spring by a factor $k_2(1-D_2)$.

After undertaking the appropriate free body analysis, with the assumptions that there is no damping and the system has no external forcing. The following equation is derived, $$\mathbf{M}\Ddot{X}+\mathbf{K}X=0$$
where,
$$
\mathbf{M}=
\begin{bmatrix}
m_1 & 0 & 0 \\
0 & m_2 & 0 \\
0 & 0 & m_3
\end{bmatrix}, \hspace{1cm}
\mathbf{K}=
\begin{bmatrix}
k_1+k_2(1-\alpha_2 T)(1-D) & -k_2(1-\alpha_2 T)(1-D) & 0 \\
-k_2(1-\alpha_2 T)(1-D) & k_2(1-\alpha_2 T)(1-D)+k_3 & -k_3 \\
0 & -k_3 & k_3 + k_4
\end{bmatrix}
$$

Assuming a harmonic response, a solution can be obtained by determining the eigenvalues of the equation $|-\mathbf{M}\omega_i^2+\mathbf{K}|=0$
The result is a matrix equation that determines the response of the system, where $\omega_i$ is a function of the parameters $\omega_i(\alpha_2,D_2,T)$. This is a $3DOF$ system, therefore, it is possible to obtain three eigenvalues $(\omega_i^2)$ and eigenvectors (mode shapes), $X_i$, where the eigenvalues and eigenvectors satisfy the equation,
$$\mathbf{M}^{-1}\mathbf{K}X=\omega^2_i X, \ i\in\{1,2,3\}$$

\subsection{TDA of Mass-Spring-Damper Model}
This example is embedded in $\mathbb{R}^4$; for this reason the topological features shown in the persistence barcodes will have to be studied to determine the topology of the data, rather than visualisation. In this case, analysis is being undertaken to determine the topology of the manifold from which the data are sampled. This was not the case with the sphere, as the topology was known prior to analysis. Despite this, it is still possible to determine some features of the topology of the data. It is known that the manifold is connected, as it is described by continuous functions. Therefore, the condition that $\beta_0=1$ must be satisfied.

It is also expected that there is a self-intersection in the manifold. This can be inferred from the model, where an increase in temperature would reduce the spring to become less stiff. This is also the case for the damage parameter. This must mean that there is a value of $\alpha T$ that will give the same reduction in $k$, that a value of $D$ will. Since all the variables are contained within $\mathbf{K}$, this means that the natural frequency must be the same value at these two points, and thus there is a self intersection in the manifold. It is actually the case that this is a family of values, and not just a single point. The family of values is parameterised by,
$$(1-\alpha_2T)(1-D_2)=(1-\alpha_2'T')(1-D_2')$$
When creating the data, the inputs were $T \in [250,500]$ with seven divisions, $\alpha_2 \in [0,0.005]$ with six divisions and $D \in [0,1]$ with six divisions; meaning that there are 252 data points. The values for $T,\alpha$ and $D$ can be thought of as coordinates of the form $(T,\alpha,D)$. Each point is then input into the equations above to give a value $\omega_i$. Each point can then be embedded in 4D in the form $(T,\alpha,D,\omega_i)$. These points will then form the vertices of the simplicial complex, which can be analysed using TDA. When generating data this way, the number of points generated is proportional to the power of the dimension, so calculations can get cumbersome for high dimensions or high sample sizes. After the points had been generated, each point was scaled down by the largest value in that dimension, meaning that all the data are inside $[0,1]^4$.

The values used for the calculations in this first example were:
$$m_1=m_2=m_3= 10kg, \ k_1=k_2=k_3=k_4=10000Nm^{-1}$$

\begin{table}[H]
    \centering
    \begin{tabular}{|c|c|c|c|}
        \hline \rowcolor[HTML]{c0c0c0}
        Description   & $\omega_1$ & $\omega_2$ & $\omega_3$ \\ \hline
        \rowcolor[HTML]{EFEFEF} 
        Epsilon      & $\epsilon_1 = 0.77$ & $\epsilon_2 = 0.33$ & $\epsilon_3 = 0.31$            \\ \hline
        Num of vert   & 252            & 252   & 252          \\ \hline
        \rowcolor[HTML]{EFEFEF} 
        Num of Simp   & 188087202          & 160639  & 93917          \\ \hline
        Betti Numbers & $[1,0,0,0]$           & $[1,0,0,0]$    & $[1,0,0,0]$         \\ \hline
        \rowcolor[HTML]{EFEFEF} 
        Dimension     & 4               &  4        & 4        \\ \hline
    \end{tabular}
        \caption{$\omega_{1,2,3}$ simplicial complexes data.}
        \label{tab:results}
\end{table}

\begin{figure}[H]
    \centering
    \includegraphics[width = \textwidth]{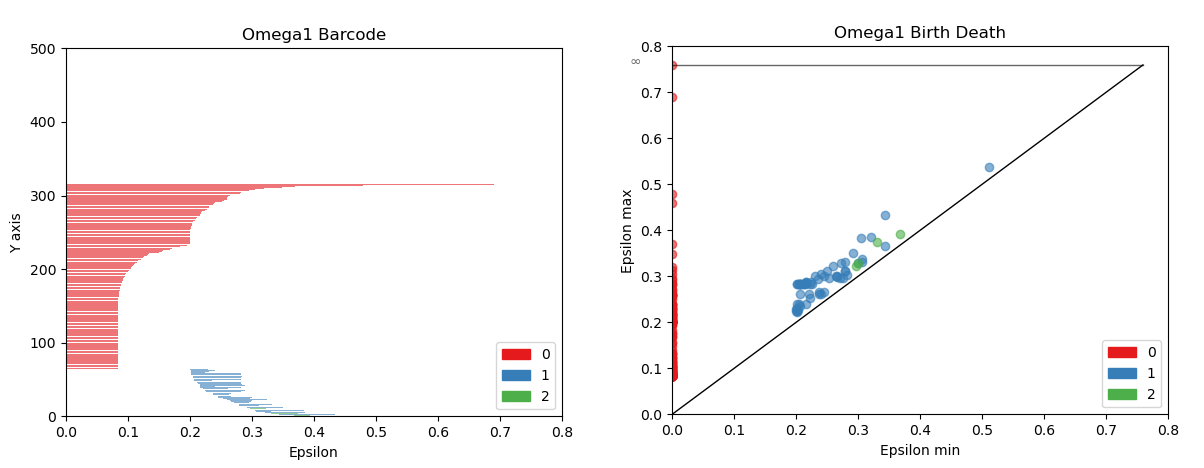}
    \caption{$\omega_1$ persistence representations.}
    \label{fig:omega1}
\end{figure}

\begin{figure}[H]
    \centering
    \includegraphics[width = \textwidth]{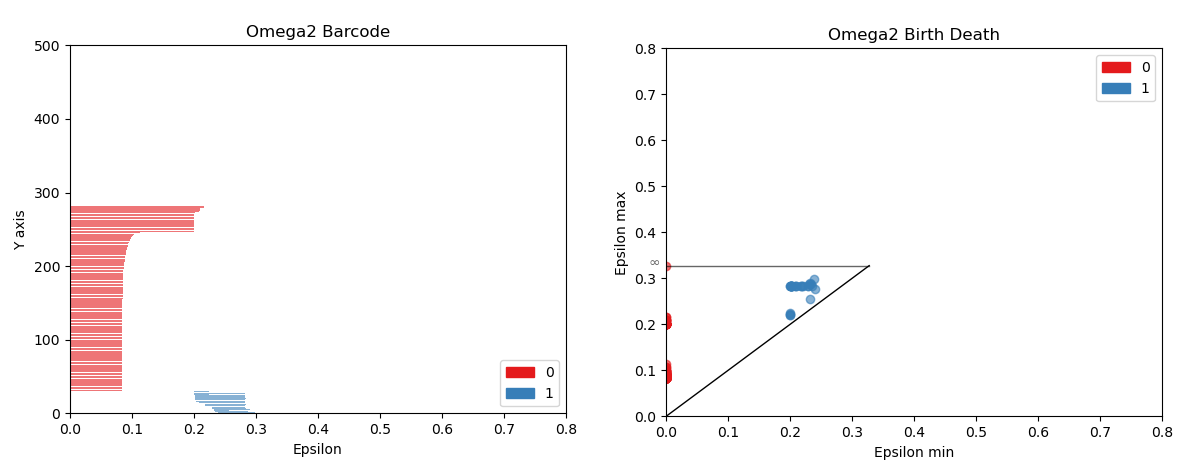}
    \caption{$\omega_2$ persistence representations.}
    \label{fig:omega2}
\end{figure}

\begin{figure}[H]
    \centering
    \includegraphics[width = \textwidth]{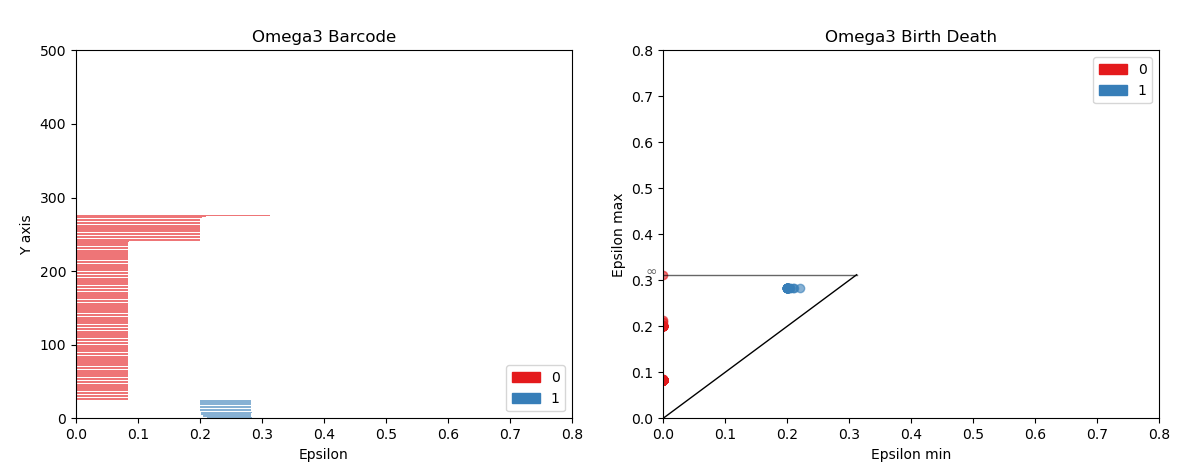}
    \caption{$\omega_3$ persistence representations.}
    \label{fig:omega3}
\end{figure}

The Betti numbers in Table \ref{tab:results} are relatively uninteresting, but these Betti numbers are for only for the values of epsilon listed. The persistence diagrams in Figures \ref{fig:omega1},\ref{fig:omega2} and \ref{fig:omega3} give more insight into the features for values of epsilon less than the ones listed in Table \ref{tab:results}.
In the case of $\omega_1$, the manifold appears to have a sudden and sharp change with the parameters; this can be seen in the birth-death diagram in Figure \ref{fig:omega1}, as the second to last red point $\approx (0,0.69)$ is the instance at which the simplicial complex is fully connected. It appears that the data are roughly divided into two clusters either side of this sudden change. A high value of epsilon is required to span this sudden change. As a result of a high value of epsilon, the clusters become highly connected within. This can be seen in the relatively large number of simplices in Table \ref{tab:results} when compared to the other natural frequencies. 

Another method of spanning the gap between the sudden change, instead of increasing $\epsilon_1$, is to increase the resolution of the data by taking more divisions for each parameter. This method also has an issue. More points are added to the sudden change, but there will also be more points in each cluster. Meaning there are more points to be connected; this massively compromises the computation time.

It can also be seen on the barcode in Figure \ref{fig:omega1} that there are numerous instances of 1D holes that persist for about 0.1 units in epsilon. As features only persist for small ranges of epsilon, this is likely caused by the lack of resolution in the data. In the case of $\omega_1$, there are also instances of cavities forming. These features also do not persist for long, and are likely artefacts due to the lack of resolution in the data.

The persistence representations shown in Figures \ref{fig:omega2} and \ref{fig:omega3} are fairly similar, both having nearly the same features present.

Both $\omega_2, \omega_3$ achieved the criteria of $\beta_0=1$ at a much lower value of $\epsilon$ when compared with $\omega_1$, which indicates there is no sudden abrupt change in the manifold. There are also 1D holes present, again this is likely to be an artefact of the lack of resolution in the data.

There is also a general trend of the persistence intervals becoming more square in the $\beta_0$ section of the barcode, going from $\omega_1$ to $\omega_3$. This says that the distance between the connected components is becoming more uniform.

\subsubsection{Values at different stiffness}
Following the previously-outlined procedure, a new data set was collected for the values,
$$m_1=m_2=m_3= 10kg, \ k_2=5000Nm^{-1}, \ k_1=k_3=k_4=10000Nm^{-1}$$
For continuity, the same number of divisions was used along each of the inputs, and the same process of selecting the lowest value of epsilon that gave a fully-connected simplicial complex. Following these conditions, the results in Table \ref{tab:results5000} were obtained.

\begin{table}[H]
    \centering
    \begin{tabular}{|c|c|c|c|}
        \hline \rowcolor[HTML]{c0c0c0}
        Description   & $\omega_1'$ & $\omega_2'$ & $\omega_3'$ \\ \hline
        \rowcolor[HTML]{EFEFEF} 
        Epsilon      & $\epsilon_1' = 0.38$ & $\epsilon_2' = 0.32$ & $\epsilon_3' = 0.30$            \\ \hline
        Num of vert   & 252            & 252   & 252          \\ \hline
        \rowcolor[HTML]{EFEFEF} 
        Num of Simp   & 339447          & 87571  & 97341          \\ \hline
        Betti Numbers & $[1,0,0,0]$           & $[1,0,0,0]$    & $[1,0,0,0]$         \\ \hline
        \rowcolor[HTML]{EFEFEF} 
        Dimension     & 4               &  4        & 4        \\ \hline
    \end{tabular}
        \caption{$\omega_{1,2,3}',\ k_2=5000Nm^{-1}$ simplicial complexes data.}
        \label{tab:results5000}
\end{table}
Immediately it is clear that the $\epsilon_1'$ is much smaller that $\epsilon_1$; this means that the computation time was much less when calculating the simplicial complex for $k_2=5000Nm^{-1}$ than $k_2=10000Nm^{-1}$.
\begin{figure}[H]
    \centering
    \includegraphics[width = 0.9\textwidth]{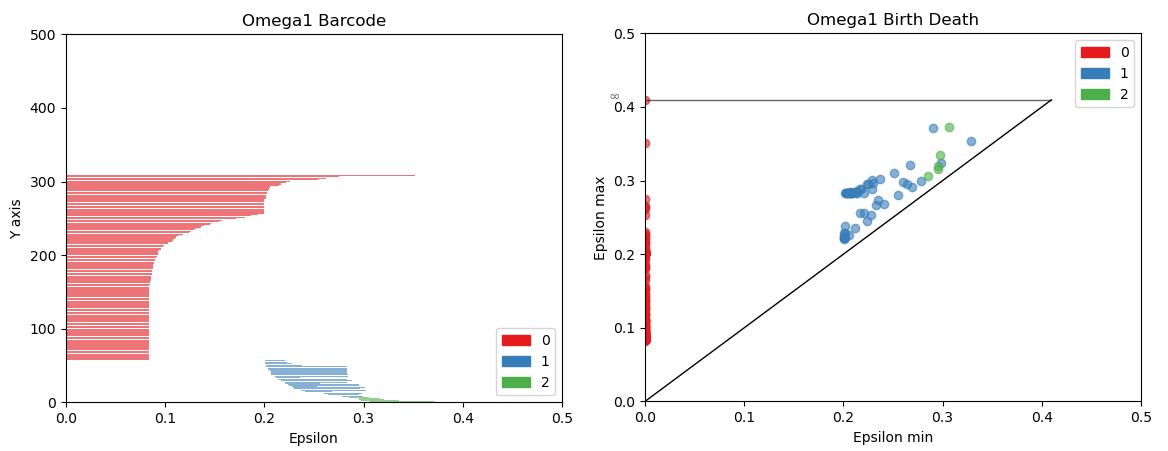}
    \caption{$\omega_1'$ persistence representations.}
    \label{fig:omega15000}
\end{figure}
\begin{figure}[H]
    \centering
    \includegraphics[width = 0.9\textwidth]{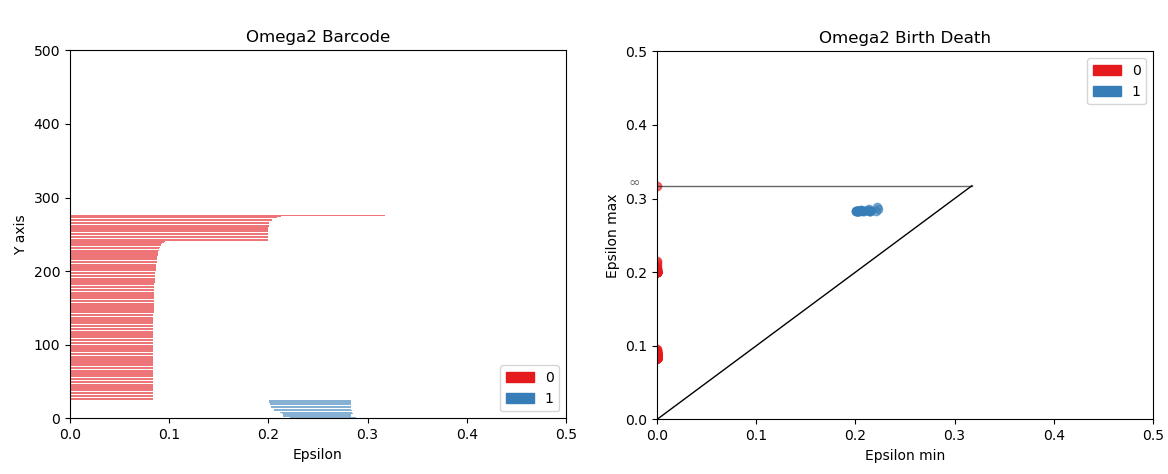}
    \caption{$\omega_2'$ persistence representations.}
    \label{fig:omega25000}
\end{figure}
\begin{figure}[H]
    \centering
    \includegraphics[width = 0.9\textwidth]{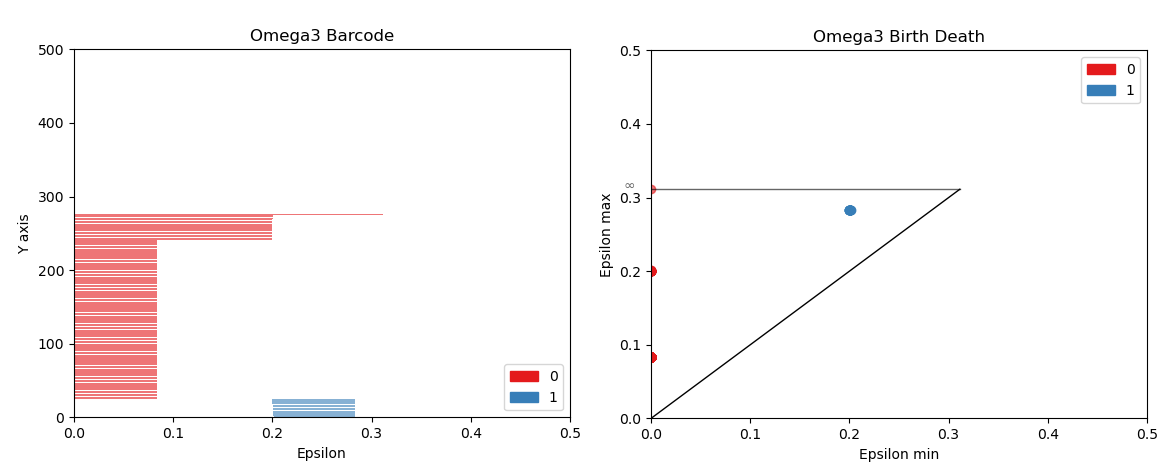}
    \caption{$\omega_3'$ persistence representations.}
    \label{fig:omega35000}
\end{figure}
The barcodes still exhibit a similar structure when compared with their counterparts in Figures \ref{fig:omega1}, \ref{fig:omega2} and \ref{fig:omega3}. When comparing the barcodes using the Wasserstein distance, the following values were obtained,
$$\partial_{W_p}(B_{\omega_1}, B_{\omega_1'}) = 0.69029,\ \partial_{W_p}(B_{\omega_2},\  B_{\omega_2'}) = 0.077248,\ \partial_{W_p}(B_{\omega_3}, B_{\omega_3'}) = 0.029308 $$

The values may appear relatively small, but the data were scaled by the largest value in that dimension, therefore, the Wasserstein distance is also scaled. There appears to be a general trend, that the higher the natural frequency, the more similar the manifolds are. From this, it can be inferred for this example, that as successive natural frequencies are computed there is a lesser reliance on the stiffness term, $k_2$.

The idea of comparing two persistence data sets is useful in engineering. If the manifold structure of data from a machine or structure operating at optimal conditions is known, the shape of the manifold could be used to fine tune or identify differences in the shape of a manifold of similar machine or structure in an unknown state. This could have interesting applications in SHM.

\section{Conclusions}
\label{sec:conclusions}
The results of this paper show that topological inference is a viable analysis strategy for engineering applications. By opting to use TDA, features in the data can be identified that would not necessarily be found in other analysis techniques. Features of the data such as holes or cavities in the underlying manifold can help to identify trends or anomalies. This inference can then be extended to determine engineering examples such as clustering data by their topological features, this can then be used to identify changes in the operating state of a machine or structure.

An appropriate method for conducting topological data analysis has been outlined, with a special consideration for conducting the analysis in an effective manner. Examples were given that highlighted where TDA can be inefficient or return misleading results.

The engineering-specific example given here is a very fundamental start, but not without its issues. Calculating data at higher dimensions is computationally expensive, and therefore the resolution in the manifold was too low. However, this is not the case with real world data that has been collected experimentally. Another issue arises when data are not evenly distributed, as this results in overly-connected areas.

In summary, a novel data analysis strategy has been discussed, previously unfamiliar to engineering, that can provide new insights into the fundamental structure and shape of the data, even in higher-dimensional analysis where an intuitive understanding of shape is no longer upheld.

\section*{Acknowledgements}
The authors would like to thank the UK EPSRC for funding via the Established Career Fellowship EP/R003645/1 and the Programme Grant EP/R006768/1.

\bibliographystyle{unsrt}
\bibliography{imac_21_TG_2}

\end{document}